\documentclass[a4,11]{article}

\makeatletter

\def\affil#1{\def\@affil{#1}}

\def\@maketitle{
\begin{center}
{\Large \@title \par}
\vspace{20pt}
{\normalsize \@author \par}
\vspace{15pt}
{\it\footnotesize \@affil \par}
\vspace{15pt}
\par\vskip 1.5em
}

\makeatother

\usepackage{amsmath,amsthm,amssymb,xypic,url}

\newtheorem{Th}{Theorem}[section]
\newtheorem{Lem}[Th]{Lemma} 
\newtheorem{Prop}[Th]{Proposition} 
\newtheorem{Cor}[Th]{Corollary}

\newtheorem{Rem}[Th]{Remark}

\def\R{{\mathbb R}}
\def\N{{\mathbb N}}
\def\Hb{{H^1(\Omega)}}
\def\H{{H^1_0(\Omega)}}
\def\Hi{{H^{-1}(\Omega)}}
\def\L{{L^2(\Omega)}}
\def\d{\operatorname{div}}
\def\into{\int_\Omega}
\def\eps{\varepsilon}
\def\DS{\displaystyle}

\title{\Large Asymptotic analysis of an $\varepsilon$-Stokes problem connecting Stokes 
and pressure-Poisson problems}

\author{\normalsize Kazunori Matsui$^{a,*}$, Adrian Muntean$^b$}
\affil{$^a$Division of Mathematical and Physical Sciences,
Graduate School of Natural Science and Technology, Kanazawa\\ 
University, Kakuma, Kanazawa 920-1192 Japan\\
$^b$Department of Mathematics and Computer Science, Karlstad University, Universitetsgatan 2, 651 88 Karlstad Sweden}

\begin{document}
\maketitle
\footnotetext[1]{Corresponding author.\\
\qquad{\it Email address:} {\tt first-lucky@stu.kanazawa-u.ac.jp} (Kazunori Matsui),
{\tt adrian.muntean@kau.se} (Adrian Muntean)\\
\qquad{\it Date:} \today}

\begin{abstract}

In this Note, we prepare an $\eps$-Stokes problem
connecting the Stokes problem and the corresponding pressure-Poisson equation
using one parameter $\eps>0$. We prove that the solution to the $\eps$-Stokes problem,
convergences as $\eps$ tends to 0 or $\infty$ to the Stokes and pressure-Poisson problem,
respectively.

\end{abstract}

%
\section{Introduction}
%

Let $\Omega$ be a bounded domain in $\R^n(n\ge 2,n\in\N)$
with Lipschitz continuous boundary $\Gamma$ and
let $F\in\L^n,u_b\in H^{1/2}(\Gamma)^n$ satisfy
$\int_\Gamma u_b\cdot \nu=0$,
where $\nu$ is the unit outward normal vector for $\Gamma$.
The weak form of the Stokes problem is : Find $u_S\in\Hb^n$ and $p_S\in\L/\R$ satisfying
\begin{align}\tag{S}
\left\{\begin{array}{ll}
-\Delta u_S+\nabla p_S = F & in~\Hi^n, \\
\d u_S = 0 & in~\L, \\
u_S=u_b & on ~H^{1/2}(\Gamma)^n.
\end{array}\right.
\end{align}
We refer to \cite{Temam} for details on the Stokes problem,
(i.e. more physical background and corresponding mathematical analysis).
Taking the divergence of the first equation, we are led to
\begin{align}\label{ppstrong}
\d F=\d(-\Delta u_S+\nabla p_S)=-\Delta (\d u_S)+\Delta p_S=\Delta p_S
\end{align}
in distributions sense.
This is often called pressure-Poisson equation and is used in MAC or SMAC method 
(cf. \cite{mac1,mac2}, e.g.).
Bearing this in mind, we consider a similar problem: 
Find $u_{PP}\in\Hb^n$ and $p_{PP}\in\Hb$ satisfying
\begin{align}\tag{PP}
\left\{\begin{array}{ll}
-\Delta u_{PP}+\nabla p_{PP} = F & in~\Hi^n, \\
-\Delta p_{PP} = -\d F& in~\Hi, \\
u_{PP}=u_b & on ~H^{1/2}(\Gamma)^n, \\
p_{PP}=p_b & on ~H^{1/2}(\Gamma).
\end{array}\right.
\end{align}
with $p_b\in H^{1/2}(\Gamma)$. Let this problem be called pressure-Poisson problem.
This idea using (\ref{ppstrong}) instead of $\d u_S=0$ is useful
to calculate the pressure numerically in the Navier-Stokes equation.
For example, the idea is used in both the MAC and SMAC methods \cite{mac1,mac2}.

\begin{figure}[h]
\[
\xymatrix{
  (PP)\ar@{.}[rr]& &(S)\\
 &(ES)\ar[ul]^{ \eps\rightarrow\infty}\ar[ur]_{\eps\rightarrow 0}&
}
\]
\caption{
  Sketch of the connections between the problems (S), (PP) and (ES).
}
\label{diagram}
\end{figure}
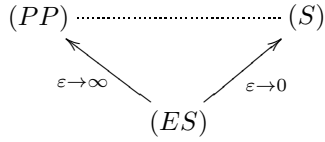

In this Note, we prepare on an ``interpolation'' between these problems (S) and (PP), 
i.e. we introduce an intermediate problem:
For $\eps>0$, find $u_\eps\in\Hb^n$ and $p_\eps\in\Hb$ which satisfy
\begin{align}\tag{ES}
\left\{\begin{array}{ll}
-\Delta u_\eps+\nabla p_\eps  = F& in~\Hi^n, \\
-\eps\Delta p_\eps+\d u_\eps = -\eps\d F& in~\Hi, \\
u_\eps=u_b & on ~H^{1/2}(\Gamma)^n, \\
p_\eps=p_b & on ~H^{1/2}(\Gamma).
\end{array}\right.
\end{align}
Let this problem be called $\eps-$Stokes problem.
The $\eps$-Stokes problem (ES) formally approximates the Stokes problem (S)
as $\eps\rightarrow0$ and the pressure-Poisson problem (PP)
as $\eps\rightarrow\infty$ (Figure \ref{diagram}).
We show here that (ES) is a natural link between (S) and (PP)
in Proposition \ref{sp_prop}.
The aim of this Note is to give a precise asymptotic estimates for (ES)
when $\eps$ tends to zero or $\infty$.

%
\section{Well-posedness}\label{sec_wellposed}
%

%
\subsection{Notation}
%
We set
\[
C^\infty_0(\Omega)^n:=\{f\in C^\infty(\Omega)^n~|~supp(f)~{\rm is~compact~subset~in~}\Omega\},
\]
\[
V:=\{u\in\H^n~|~\d u=0\},~
\L/\R:=\{u\in\L~|~\into u=0\}.
\]

For $m=1$ or $n$, $\Hi^m=(\H^m)^*$ is equipped with the norm
$
||f||_{\Hi^m}:=\sup_{\varphi\in S_m}\langle f,\varphi\rangle
$
for $f\in\Hi^m$,
where $S_m=\{\varphi\in\H^m~|~||\nabla\varphi||_{\L^{n\times m}}=1\}$.
We define $[p]:=p-(1/|\Omega|)\into p$ and
$||p||_{\L/\R}:=\inf_{a\in\R}||p-a||_\L=||[p]||_\L$ for all $p\in\L$,
where $|\Omega|$ is the volume of $\Omega$.

Let $\gamma_0\in B(\Hb, H^{1/2}(\Gamma))$ be
the standard trace operator.
It is known that (see e.g. \cite[p.10,11,Lemma 1.3]{Temam}) 
there exists a linear continuous operator 
$\gamma_\nu :\Hb^n\rightarrow H^{-1/2}(\Gamma)$ such that
$\gamma_\nu u=u\cdot \nu |_\Gamma$ for all $u\in C^\infty(\overline{\Omega})^n$,
where $\nu$ is the unit outward normal for $\Gamma$ 
and $H^{-1/2}(\Gamma):=H^{1/2}(\Gamma)^*$.  
Then, the following generalized Gauss divergence formula holds:
\[
\into u\cdot\nabla \omega+\into(\d u)\omega=\langle\gamma_\nu u,\gamma_0 \omega\rangle
\quad{ for~all~}u\in \Hb^n,\omega\in\Hb.
\]

We recall the following Theorem \ref{lem_Necas}
that plays an important role in the proof of 
the existence of pressure solution of Stokes problem;
see \cite[p.187-190, Lemme 7.1, $l=0$]{Necas} and 
\cite[p.111-115, Theorem 3.2 and Remark 3.1 ($\Omega$ is $C^1$ class)]{Duvaut} for the proof.
\begin{Th}\label{lem_Necas}
  There exists a constant $c>0$ such that
  \[
  ||f||_\L\le c(||f||_\Hi+||\nabla f||_\Hi)
  \]
  for all $f\in\L$.
\end{Th}

The following result follows from Theorem \ref{lem_Necas}.

\begin{Th}\label{grad}
  {\rm \cite[p.20-21]{Girault}}
  There exists a constant $c>0$ such that
  \[
  ||f||_{\L/\R}\le c||\nabla f||_{\Hi^n}
  \]
  for all $f\in\L$.
\end{Th}

%
\subsection{Well-posedness}
%

\begin{Th}\label{stokes_t}
For $F\in\L^n$ and $u_b\in H^{1/2}(\Gamma)^n$,
there exists a unique pair of functions $(u_S,p_S)\in \Hb^n\times(\L/\R)$ satisfying (S).
\end{Th}

See \cite[p.31-32,Theorem 2.4 and Remark 2.5]{Temam} for the proof.

\begin{Th}\label{pp_thm}
For $F\in\L^n,u_b\in H^{1/2}(\Gamma)^n$ and $p_b\in H^{1/2}(\Gamma)$,
there exists a unique pair of functions $(u_{PP},p_{PP})\in \Hb^n\times\Hb$ satisfying (PP).
\end{Th}

\noindent ${\textit Proof.}$
From the second and fourth equations of (PP), $p_{PP}\in\Hb$ is uniquely determined.
Then $u_{PP}\in\Hb^n$ is also uniquely determined from the first and third equations.
\qed

\begin{Cor}\label{sp_cor}
If the solution $(u_{PP},p_{PP})\in \Hb^n\times\Hb$ of (PP) satisfies
$\d u_{PP}=0$, by Theorem \ref{stokes_t}, $u_S=u_{PP}$ and $p_S=[p_{PP}]$ hold.
\end{Cor}

\begin{Th}\label{estokes_thm}
For $\eps>0,F\in\L^n,u_b\in H^{1/2}(\Gamma)^n$ and $p_b\in H^{1/2}(\Gamma)$,
there exists a unique pair of functions $(u_\eps,p_\eps)\in \Hb^n\times\Hb$ 
satisfying the problem (ES).
\end{Th}

\noindent ${\textit Proof.}$
We pick $u_1\in\Hb^n$ and $p_0\in\Hb$ with $\gamma_0 u_1=u_b,\gamma_0 p_0=p_b$.
Since $\d:\H^n\rightarrow\L/\R$ is surjective
\cite[p.24, Corollary 2.4, 2${}^\circ$)]{Girault},
there exists $u_2\in\H^n$ such that
$\d u_2=\d u_1$. We put $u_0:=u_1-u_2$, and then $\gamma_0 u_0=u_b$ and
$\d u_0=0$ in $\Omega$.
To simplify the notation, we set
$u:=u_\eps-u_0(\in\H^n),p:=p_\eps-p_0(\in\H),f\in\Hi^n$ and $g\in\Hi$ such that
$\langle f,v\rangle=\into Fv-\into\nabla u_0:\nabla v-\into(\nabla p_0)\cdot v
~(v\in\H^n),
\langle g,q\rangle=\into F\cdot\nabla q-\into\nabla p_0\cdot\nabla q~(q\in\H)$.
Then we have 
\begin{align}\label{estokes_exist}
  \left\{\begin{array}{ll}
    {\DS
      \into\nabla u:\nabla \varphi +\into (\nabla p)\cdot\varphi
      =\langle f,\varphi\rangle  }
    & {\rm for~all~} \varphi\in \H^n,\\
    {\DS
      \eps\into\nabla p\cdot\nabla\psi +\into(\d u)\psi
      =\eps\langle g,\psi\rangle }
    & {\rm for~all~} \psi \in\H.
  \end{array}\right.
\end{align}
Adding the equations in (\ref{estokes_exist}), we get
\[
(\left( \begin{array}{c}u\\ p\end{array}\right) ,
  \left( \begin{array}{c}\varphi \\ \psi\end{array}\right) )_\eps
    = \langle f,\varphi\rangle +\eps\langle g,\psi\rangle.
\]
Here, we denote
\[
(\left( \begin{array}{c}u\\ p\end{array}\right) ,
 \left( \begin{array}{c}\varphi \\ \psi\end{array}\right) )_\eps
   :=\into\nabla u:\nabla\varphi + \eps\into \nabla p\cdot \nabla \psi
   +\into(\nabla p)\cdot\varphi + \into (\d u)\psi.
\]
We check that $(*,*)_\eps$ is a continuous coercive bilinear form
on $\H^n\times\H$. The bilinearity and continuity of $(*,*)_\eps$ are obvious.
The coercivity of $(*,*)_\eps$ is obtained in the following way:
Let ${}^t(u,p)\in\H^n\times\H$. We have the following sequence of inequalities;
\begin{align*}
  \begin{array}{rl}
    (\left( \begin{array}{c}u\\ p\end{array}\right) ,
      \left( \begin{array}{c}u\\ p\end{array}\right) )_\eps
        =& {\DS \into\nabla u:\nabla u + \eps\into \nabla p\cdot \nabla p
        +\into \d  (up)}\\
        =& ||\nabla u||^2_{\L}+\eps||\nabla p||^2_{\L}\\
        \ge& {\rm min}\{ 1,\eps\}(||\nabla u||^2_{\L}+||\nabla p||^2_{\L})\\
        \ge& c~{\rm min}\{ 1,\eps\}(||u||^2_{\Hb^n}+||p||^2_{\Hb})
  \end{array}
\end{align*}
by the Poincar\' e inequality.
Summarizing, $(*,*)_\eps$ is a continuous coercive bilinear form and
$H^1_0(\Omega)^{n+1}$ is a Hilbert space.
Therefore, the conclusion of Theorem \ref{estokes_thm} follows based on
the Lax-Milgram Theorem.
\qed

From now on, let the solutions of (S), (PP) and (ES) be denoted by
$(u_S,p_S),(u_{PP},p_{PP})$ and $(u_\eps,p_\eps)$, respectively.

\begin{Prop}\label{sp_prop}
Suppose that $p_S\in\Hb$.
Then there exists a constant $c>0$ independent of $\eps$ such that
\[
||u_S-u_{PP}||_{\Hb^n}\le c||\gamma_0 p_S-p_b||_{H^{1/2}(\Gamma)},\quad
||u_S-u_\eps||_{\Hb^n}\le c||\gamma_0 p_S-p_b||_{H^{1/2}(\Gamma)}.
\]
In particular, if $\gamma_0p_S=p_b$,
then $p_{PP}=p_\eps=p_S$ hold for all $\eps>0$.
\end{Prop}

\noindent ${\textit Proof.}$
From (S) and (PP), we have 
\begin{align}\label{s-p} \left\{\begin{array}{ll}{\DS
\into\nabla (u_S-u_{PP}):\nabla\varphi 
=-\into (\nabla(p_S-p_{PP}))\cdot\varphi }
& {\rm for~all~} \varphi \in \H^n, \\
{\DS
\into\nabla(p_S-p_{PP})\cdot\nabla\psi=0}
&{\rm for~all~} \psi \in \H.
\end{array}\right.\end{align}
Putting $\varphi:=u_S-u_{PP}\in\H^n$ in (\ref{s-p}), we get
\[\begin{array}{rl}
||\nabla(u_S-u_{PP})||^2_{\L^{n\times n}}
&={\DS -\into (\nabla(p_S-p_{PP}))\cdot(u_S-u_{PP})}\\
&\le ||\nabla(p_S-p_{PP})||_{\L^n}||u_S-u_{PP}||_{\L^n},
\end{array}\]
and then $||u_S-u_{PP}||_{\Hb^n}\le c_1||\nabla(p_S-p_{PP})||_{\L^n}$ follows.
From the second equation of (\ref{s-p}),
there exists a constant $c_2>0$ such that
$||\nabla(p_S-p_{PP})||_{\L^n}\le c_2||\gamma_0 p_S-p_b||_{H^{1/2}(\Gamma)}$.
Therefore we obtain
$
||u_S-u_{PP}||_{\Hb^n}\le c_1c_2||\gamma_0 p_S-p_b||_{H^{1/2}(\Gamma)}.
$

Let $w_\eps:=u_S-u_\eps\in\H^n,r_\eps:=p_{PP}-p_\eps\in\H$.
By (S), (PP) and (ES), we have
\begin{align}\label{s-p2} \left\{\begin{array}{ll}{\DS
    \into\nabla w_\eps:\nabla\varphi 
    +\into(\nabla r_\eps)\cdot\varphi 
    =-\into(\nabla(p_S-p_{PP}))\cdot\varphi }
  & {\rm for~all~} \varphi \in \H^n, \\
  {\DS
    \eps\into\nabla r_\eps\cdot\nabla\psi+\into(\d w_\eps)\psi=0}
  &{\rm for~all~} \psi \in \H.
  \end{array}\right.\end{align}
Putting $\varphi:=w_\eps$ and $\psi:=r_\eps$ and adding two equations of (\ref{s-p2}), 
we get
\[
||\nabla w_\eps||^2_{\L^{n\times n}}+\eps||\nabla r_\eps||^2_{\L^n}
=-\into(\nabla(p_S-p_{PP}))\cdot w_\eps
\le ||\nabla(p_S-p_{PP})||_{\L^n}||w_\eps||_{\L^n}
\]
from $\into(\nabla r_\eps)\cdot w_\eps=-\into(\d w_\eps)r_\eps$.
Thus it leads $||w_\eps||_{\Hb^n}\le c_3||\nabla(p_S-p_{PP})||_{\L^n}$.
Hence we obtain
$
||u_S-u_\eps||_{\Hb^n}=||w_\eps||_{\Hb^n}\le c_2c_3||\gamma_0 p_S-p_b||_{H^{1/2}(\Gamma)}.
$
\qed

\begin{Prop}\label{sp_prop2}
  Under the hypotheses of Proposition \ref{sp_prop},
  if $\tilde p\in\Hb$ satisfies $\gamma_0\tilde p=p_b$, then we have
  \[
  ||\nabla(\tilde p-p_{PP})||_{\L^n}\le||\nabla(\tilde p-p_S)||_{\L^n}.
  \]
\end{Prop}

\noindent ${\textit Proof.}$
By the second equation of (\ref{s-p}) and $\tilde p-p_{PP}\in\H$, it yields
\[\begin{array}{rl}
||\nabla(\tilde p-p_{PP})||^2_{\L^n}
&={\DS \into\nabla(\tilde p-p_S+p_S-p_{PP})\cdot\nabla(\tilde p-p_{PP})}\\
&\le||\nabla(\tilde p-p_S)||_{\L^n}||\nabla(\tilde p-p_{PP})||_{\L^n}.\\
\end{array}\]
Hence we obtain $||\nabla(\tilde p-p_{PP})||_{\L^n}\le||\nabla(\tilde p-p_S)||_{\L^n}$.
\qed

\begin{Rem}
If $p_S\in\Hb$, then we have
\[
||\nabla(p_\eps-p_{PP})||_{\L^n}\le||\nabla(p_\eps-p_S)||_{\L^n}
\]
for all $\eps>0$, (from Proposition \ref{sp_prop2}). 
Hence, if $(\nabla p_\eps)_{\eps>0}$ converges strongly to $\nabla p_S$ in $\L^n$, 
then $u_{PP}=u_S$ and $p_{PP}=p_S$, which imply $\gamma_0 p_S=p_b$.
In other words, if $p_S\in\Hb$ satisfies $\gamma_0 p_S\ne p_b$,
then $\nabla p_\eps$ does not converge to $\nabla p_S$ in $\L^n$ as $\eps\rightarrow0$.
\end{Rem}

%
\section{Links between (ES) and (PP)}\label{sec_ep}
%

\begin{Th}\label{ep_conv}
There exists a constant $c>0$ independent of $\eps>0$ satisfying
\[
||u_\eps-u_{PP}||_{\Hb^n}\le\frac{c}{\eps}||\d u_{PP}||_\Hi,\quad
||p_\eps-p_{PP}||_\Hb\le\frac{c}{\eps}||\d u_{PP}||_\Hi.
\]
for all $\eps>0$. In particular, we have
\[
||u_\eps-u_{PP}||_{\Hb^n}\rightarrow 0,~
||p_\eps-p_{PP}||_\Hb\rightarrow 0
\quad as~\eps\rightarrow\infty.
\]
\end{Th}

\noindent ${\textit Proof.}$
From (PP) and (ES), we have
\begin{align}\label{esp}
  \left\{\begin{array}{ll}
    {\DS
      \into\nabla (u_\eps-u_{PP}):\nabla\varphi
      +\into (\nabla(p_\eps-p_{PP}))\cdot\varphi =0}
    & {\rm for~all~} \varphi\in\H^n,\\
    {\DS
      \eps\into\nabla(p_\eps-p_{PP})\cdot\nabla\psi+\into(\d (u_\eps-u_{PP}))\psi
	=-\into(\d u_{PP})\psi}
    & {\rm for~all~} \psi \in\H.
  \end{array}\right.
\end{align}
Putting $\varphi:=u_\eps-u_{PP}\in\H^n$ and $\psi:=p_\eps-p_{PP}\in\H$ and
adding two equations of (\ref{esp}), we obtain
\[
||\nabla (u_\eps-u_{PP})||^2_{\L^{n\times n}}
+\eps||\nabla (p_\eps-p_{PP})||^2_{\L^n}
\le ||\d u_{PP}||_\Hi||\nabla (p_\eps-p_{PP})||_{\L^n},
\]
where we have used
$\into(\nabla(p_\eps-p_{PP}))\cdot(u_\eps-u_{PP})
=-\into(\d(u_\eps-u_{PP}))(p_\eps-p_{PP}).$
Thus 
\[
||\nabla (p_\eps-p_{PP})||_{\L^n}\le \frac{1}{\eps}||\d u_{PP}||_\Hi
\]
follows. In addition, by (\ref{esp}) and the Poincar\'e inequality, we have 
\begin{align*}
  \begin{array}{rl}
  ||\nabla (u_\eps-u_{PP})||^2_{\L^n}
  = &-\into(\nabla (p_\eps-p_{PP}))\cdot (u_\eps-u_{PP})\\
  \le &||\nabla (p_\eps-p_{PP})||_{\L^n}||u_\eps-u_{PP}||_{\L^n}\\
  \le &c||\nabla (p_\eps-p_{PP})||_{\L^n}||\nabla (u_\eps-u_{PP})||_{\L^n},
  \end{array}
\end{align*}
and then
$
||\nabla (u_\eps-u_{PP})||_{\L^n}\le (c/\eps)||\d u_{PP}||_\Hi
$
follows.
\qed

\begin{Cor}
If $u_{PP}$ satisfies $\d u_{PP}=0$,
then $u_\eps=u_{PP}$ and $p_\eps=p_{PP}$
hold for all $\eps>0$. Furthermore, 
$u_S=u_\eps=u_{PP}$ and $p_S=[p_\eps]=[p_{PP}]$ hold for all $\eps>0$.
\end{Cor}

%
\section{Links between (ES) and (S)}\label{sec_es}
%

\begin{Lem}\label{pleu}
  If $v\in\Hb^n,q\in\L$ and $f\in\Hi^n$ satisfy
  \[
  \into\nabla v:\nabla\varphi +\langle\nabla q,\varphi\rangle=\langle f,\varphi\rangle
  \quad for~all~\varphi\in\H^n,
  \]
  then there exists a constant $c>0$ such that
  \[
  ||q||_{\L/\R}\le c(||\nabla v||_{\L^{n\times n}}+||f||_{\Hi^n}).
  \]
\end{Lem}

\noindent ${\textit Proof.}$
Let $c$ be the constant arising in Theorem \ref{grad}. Then we have
\begin{align*}
  \begin{array}{rl}
  ||q||_{\L/\R}
  \le&c||\nabla q||_{\Hi^n}
  ={\DS c\sup_{\varphi\in S_n}
  |\langle\nabla q,\varphi\rangle|}\\
  \le&{\DS c\sup_{\varphi\in S_n}
  (|\into \nabla v:\nabla\varphi|+|\langle f,\varphi\rangle|)}\\
  \le&c(||\nabla v||_{\L^{n\times n}}+||f||_{\Hi^n}).
  \end{array}
\end{align*}
\qed

\begin{Th}\label{bb}
There exists a constant $c>0$ independent of $\eps$ such that
\[
||u_\eps||_{\Hb^n}\le c,\quad||p_\eps||_{\L/\R}\le c
\quad{for~all~\eps>0}.
\]
Furthermore, we have
\[
u_{\eps}-u_S\rightharpoonup 0~weakly~in~\H^n,~
[p_{\eps}]-p_S\rightharpoonup 0~weakly~in~\L/\R
\quad{\rm as}~\eps\rightarrow 0.
\]
\end{Th}

\noindent ${\textit Proof.}$
We use the notations $u_0\in\Hb^n,p_0\in\Hb,f\in\Hi^n$ and $g\in\Hi$ 
in Theorem \ref{estokes_thm}. We put
$\tilde u_\eps:=u_\eps-u_0\in\H^n,\tilde p_\eps:=p_\eps-p_0\in\H$.
Then we have
\begin{align}\label{es_hform}
  \left\{\begin{array}{ll}
    {\DS
      \into\nabla\tilde u_\eps:\nabla \varphi +\into (\nabla\tilde p_\eps)\cdot\varphi
      =\langle f,\varphi\rangle  }
    & {\rm for~all~} \varphi\in \H^n,\\
    {\DS
      \eps\into\nabla\tilde p_\eps\cdot\nabla\psi +\into(\d \tilde u_\eps)\psi
      =\eps\langle g,\psi\rangle }
    & {\rm for~all~} \psi \in\H.
  \end{array}\right.
\end{align}
Putting $\varphi:=\tilde u_\eps,\psi:=\tilde p_\eps$ and adding the two equations of
(\ref{es_hform}), we get
\[
||\nabla\tilde u_\eps||^2_{\L^{n\times n}}+\eps||\nabla\tilde p_\eps||^2_{\L^n}
\le ||f||_{\Hi^n}||\nabla\tilde u_\eps||_{\L^{n\times n}}
+\eps ||g||_\Hi||\nabla\tilde p_\eps||_{\L^n}
\]
since $\into(\nabla\tilde p_\eps)\cdot\tilde u_\eps=-\into(\d\tilde u_\eps)\tilde p_\eps$.
It leads that $(||\tilde u_\eps||_{\Hb^n})_{0<\eps<1}$ and
$(||\sqrt{\eps}\tilde p_\eps||_\Hb)_{0<\eps<1}$ are bounded. In addition,
\[
||\tilde p_\eps||_{\L/\R}\le{c}(||\nabla\tilde u_\eps||_{\L^{n\times n}}+||f||_{\Hi^n})
\]
by Lemma \ref{pleu}, which implies that $(||\tilde p_\eps||_{\L/\R})_{0<\eps<1}$ is bounded.
By Theorem \ref{ep_conv}, $(||u_\eps||_{\Hb^n})_{\eps\ge 1}$ and
$(||\tilde p_\eps||_{\L/\R})_{\eps\ge 1}$ are bounded, and then
$(||u_\eps||_{\Hb^n})_{\eps>0}$ and $(||\tilde p_\eps||_{\L/\R})_{\eps>0}$
are bounded.

Since $\H^n\times\L/\R$ is reflexive and 
$(\tilde u_\eps,[\tilde p_\eps])_{0<\eps<1}$ is bounded in $\H^n\times\L/\R$, 
there exist $(u,p)\in\H^n\times(\L/\R)$ and
a subsequence of pair $(\tilde u_{\eps_k}, \tilde p_{\eps_k})_{k\in\N}\subset \H^n\times\H$ 
such that
\[
\tilde u_{\eps_k}\rightharpoonup u~{\rm weakly~in~}\H^n,~
[\tilde p_{\eps_k}]\rightharpoonup p~{\rm weakly~in~}\L/\R\quad{\rm as~}k\rightarrow\infty.
\]
Hence, from (\ref{es_hform}) with $\eps:=\eps_k$, taking $k\rightarrow\infty$,
we obtain
\begin{align}\label{eslimit}
\left\{ \begin{array}{ll}{\displaystyle
\into\nabla u:\nabla\varphi +\langle \nabla p,\varphi\rangle
=\langle f,\varphi\rangle }
& {\rm for~all~} \varphi\in \H^n\\
{\displaystyle \into(\d u)\psi=0}& {\rm for~all~} \psi \in\H,
\end{array}\right.
\end{align}
where
\[
|\eps_k\into\nabla\tilde p_{\eps_k}\cdot\nabla\psi|
\le \sqrt{\eps_k}||\sqrt{\eps}\tilde p_\eps||_\Hb||\psi||_\Hb\rightarrow 0,
\]
\[
\into\nabla\tilde p_{\eps_k}\cdot\varphi
=-\into[\tilde p_{\eps_k}]\d\varphi
\rightarrow-\into p\d\varphi=\langle\nabla p,\varphi\rangle
\]
as $k\rightarrow\infty$.
The first equation of (\ref{eslimit}) implies that
\[
-\Delta(u+u_0)+\nabla(p+p_0)=F\quad{\rm in~}\Hi^n.
\]
From the second equation of (\ref{eslimit}), $\d(u+u_0)=0$ follows.
Hence, we obtain $u_S=u+u_0$ and $p_S=p+[p_0]$. Then we have 
\[
u_\eps-u_S=u_\eps-u-u_0=\tilde u_{\eps_k}-u_S\rightharpoonup 0{\rm~weakly~in~}\H^n,
\]
\[
[p_\eps]-p_S=[p_\eps-p-p_0]=[\tilde p_{\eps_k}]-p\rightharpoonup 0{\rm~weakly~in~}\L/\R.
\]
An arbitrarily chosen subsequence of $((u_\eps,[p_\eps]))_{0<\eps<1}$
has a subsequence which converges to $(u_S,p_S)$,
so we can conclude the proof.
\qed

\begin{Th}
We have
\[
u_{\eps}-u_S\rightarrow 0~strongly~in~\H^n,~
[p_{\eps}]-p_S\rightarrow 0~strongly~in~\L/\R
\quad{\rm as}~\eps\rightarrow 0.
\]
\end{Th}

\noindent ${\textit Proof.}$
We have
\begin{align}\label{ess}
  \left\{
  \begin{array}{ll}
    {\displaystyle
  \into\nabla (u_\eps-u_S):\nabla\varphi +\into (\nabla(p_\eps-p_S))\cdot\varphi =0  }
    & {\rm for~all~} \varphi\in\H^n,\\
    {\displaystyle
    \eps\into\nabla(p_\eps-p_S)\cdot\nabla\psi+\into(\d u_\eps)\psi=0}
    & {\rm for~all~} \psi \in\H.
  \end{array}
  \right.
\end{align}
We use the notations $p_0\in\Hb$ in Theorem \ref{estokes_thm}.
Putting $\varphi:=u_\eps-u_S\in\H^n,\psi:=p_\eps-p_0\in\H$ and
$\tilde p_S:=p_S-p_0\in\Hb$, we get
\[\begin{array}{ll}
&||\nabla(u_\eps-u_S)||^2_{\L^{n\times n}}+\eps||\nabla (p_\eps-p_S)||^2_{\L^n}\\
=&{\DS \into (\nabla\tilde p_S)\cdot (u_\eps-u_S)-\eps\into \nabla (p_\eps-p_S)\cdot \nabla\tilde p_S}\\
\le &||\nabla\tilde p_S||_{\L^n}||u_\eps-u_S||_{\L^n}
+\eps||\nabla (p_\eps-p_S)||_{\L^n}||\nabla\tilde p_S||_{\L^n}.\\
\end{array}\]
By Corollary \ref{bb} and the Rellich-Kondrachov Theorem, 
there exists a sequence $(\eps_k)_{k\in\N}\subset \R$
such that
\[
u_{\eps_k}\rightarrow u_S~strongly~in~\L^n.
\]
So, we can write that
\[
||\nabla(u_{\eps_k}-u_S)||^2_{\L^{n\times n}}
\le||\nabla\tilde p_S||_{\L^n}||u_{\eps_k}-u_S||_{\L^n}
+\sqrt{\eps_k}||\sqrt{\eps_k}\nabla (p_{\eps_k}-p_S)||_{\L^n}||\nabla\tilde p_S||_{\L^n}
\longrightarrow 0.
\]
It implies that
\[
||[p_{\eps_k}]-p_S||_{\L/\R}=||p_{\eps_k}-p_S||_{\L/\R}\le c||\nabla(u_{\eps_k}-u_S)||_{\L^{n\times n}}
\longrightarrow 0
\]
by Lemma \ref{pleu}.
An arbitrarily chosen subsequence of $((u_\eps,[p_\eps]))_{0<\eps<1}$
has a subsequence which converges to $(u_S,p_S)$,
so we can conclude the proof.
\qed

\bibliography{estokes_CRM}

\begin{thebibliography}{1}
\expandafter\ifx\csname url\endcsname\relax
  \def\url#1{\texttt{#1}}\fi
\expandafter\ifx\csname urlprefix\endcsname\relax\def\urlprefix{URL }\fi
\expandafter\ifx\csname href\endcsname\relax
  \def\href#1#2{#2} \def\path#1{#1}\fi

\bibitem{Temam}
R.~Temam, Navier-Stokes Equations, North Holland, 1979.

\bibitem{mac1}
F.~H. Harlow, J.~E. Welch, Numerical calculation of time-dependent viscous
  incompressible flow of fluid with a free surface, The Physics of Fluids 8
  (1965) 2182--2189.

\bibitem{mac2}
S.~McKee, M.~F. Tom\'e, J.~A. Cuminato, A.~Castelo, V.~G. Ferreira, Recent
  advances in the marker and cell method, Arch. Comput. Meth. Engng. 2 (2004)
  107--142.

\bibitem{Necas}
J.~Ne{\u{c}}as, Les m\'ethodes directes en th\'eorie des \'equations
  elliptiques, Academia, Praha, and Masson et Cie, Editeurs, Paris, 1967.

\bibitem{Duvaut}
G.~Duvaut, J.~L. Lions, Inequalities in Mechanics and Physics, Springer-Verlag,
  Berlin, 1976.

\bibitem{Girault}
V.~Girault, P.-A. Raviart, Finite Element Methods for Navier-Stokes Equations,
  Springer-Verlag, 1986.

\end{thebibliography}
\end{document}